

\magnification=\magstep1
\baselineskip=16 pt
\mathsurround=2 pt


\def\gr{\it}
\def\tit#1{\bigbreak\noindent{\bf #1}\medskip\nobreak}
\def\(#1){\smallskip\itemitem{(#1)}}
\def\PF{\noindent{\it Proof:\/} }
\def\FP{\hbox{\quad\vrule height6pt width4pt depth0pt}}

\def\IFF{if and only if}
\def\AND{{\rm\ and\ }}
\def\OR{{\rm\ or\ }}
\def\IMP{\Rightarrow}
\def\V{\,\forall\,}
\def\E{\,\exists\,}

\def\<#1>{\langle#1\rangle}
\def\EXT#1{\,\overline{#1}\,}
\def\K#1{{\rm K}(#1)}
\def\sc{\subset}

\def\R{{\bf R}}
\def\N{{\bf N}}
\def\Z{{\bf Z}}

\def\A{{\rm A}}
\def\B{{\rm B}}
\def\C{{\rm C}}
\def\D{{\rm D}}
\def\I{\rm I}

\def\:{\colon}
\def\o{\circ}
\def\P{{\rm P}}
\def\f{f}
\def\g{g}

\def\RR{\EXT\R}
\def\ZZ{\EXT\Z}

\def\AA{\EXT\A}
\def\BB{\EXT\B}
\def\CC{\EXT\C}
\def\DD{\EXT\D}
\def\II{\EXT\I}
\def\PP{\EXT\P}
\def\ff{\EXT\f}

\def\a{\alpha}
\def\b{\beta}
\def\ga{\gamma}
\def\dt{\delta}
\def\ep{\varepsilon}
\def\l{\lambda}
\def\m{\mu}
\def\n{\nu}
\def\x{\xi}
\def\r{\rho}
\def\s{\sigma}
\def\y{\upsilon}
\def\ps{\psi}
\def\w{\omega}

\def\8{\infty}
\def\d{\partial}

\def\={\approx}
\def\@{\sim}
\def\OO#1{{\cal O}(#1)}


\hrule height0pt
\vfil
{\bf VIRTUAL CALCULUS --- PART I}
\bigskip\bigskip\bigskip
{\bf S\'ergio F.~Cortizo}
\bigskip\bigskip
Instituto de Matem\'atica e Estat\'\i stica, Universidade de S\~ao
Paulo

Cidade Universit\'aria, Rua do Mat\~ao, 1010

05508-900, S\~ao Paulo, SP, Brasil
\smallskip
cortizo@ime.usp.br
\bigskip\bigskip\bigskip
{\bf Abstract}
\bigskip
\item{}An alternative organization for Differential and Integral Calculus,
based on an extension of real numbers that include infinitesimal and infinite
quantities, is presented. Only Elementary Set Theory is used, without
reference to methods or results from Mathematical Logic.
\smallskip
PACS \ 02.90.+p
\vfil\eject

\tit{I. Introduction}

The title of this work refers to an alternative organization for Differential
and Integral Calculus, based on an extension of the real numbers set, not on
the traditional concept of ``limit''. We discuss the relevancy of such
alternative below.

The notion of limit, traditionally defined by manipulation of logical
quantifiers, has the historical merit for having risen the foundations of
Calculus to the level of formal rigour required by contemporary mathematics.
However, the use of that notion as a base of Calculus can be criticized for
many reasons:

(i) The formal definition of limit is ``non-intuitive'', i.e., there is a
large distance between the rigorous definition and the corresponding
heuristic concept. The use of a suggestive notation, which introduces a
``notion of time'' on the strictly symbolic plane, does not seem to solve
this problem satisfactorily.

(ii) The formal demonstrations of results involving the concept of limit are
clumsy and artificial, due to the simultaneous use of many logic quantifiers
in a same proposition. We know that ``technicalities'' are sometimes
unavoidable in formal mathematics. Nevertheless, the amount present in the
manipulation of limits is not reasonable, obscuring even the proof of the
most elementary results.

(iii) The notation and operational rules about limits are not comfortable for
``calculations'', and they exclude many important notions which are kept
undefined from the formal viewpoint. Perhaps the best example is the concept
of ``order of an infinitesimal'' which is essential for comprehension of many
important mathematical facts.

The main motivation for this work was to approximate the intuitive ideas of
Calculus to its formalizations. Having done so, we expect that definitions
and statements become simpler and clearer, without loss of rigour which
characterizes contemporary mathematical methods. Moreover, we hope to
approximate formal proofs to the intuitive ideas behind the results.

To accomplish that, we will apply the process of ``virtual extension''
described in Ref.~1 to the real numbers set. We do not intend to show here a
complete and detailed formulation of Calculus based on this extension, since
it would be too long. Our aim is just to indicate how it can be done, and
offer enough elements so that this Calculus alternative organization can be
compared to the traditional one.
\vfill\eject
\tit{II. Virtual Numbers}

We will call {\gr virtual numbers}, or simply {\gr virtuals}, the elements of
the virtual extension $\RR$ of the set~$\R$ of ``real numbers''. According to
the identification $\R=\K\R$, we will consider $\R\sc\RR$, so that every real
number is a virtual number. (The reciprocal statement clearly does not hold.)

The real numbers will be represented by low-case Latin letters ($x$, $y$,
$a$, $b\ldots$), and generic virtual numbers by low-case Greek letters ($\x$,
$\a$, $\b$, $\l$, $\w\ldots$). The letter `$\pi$' is an exception, and it
will keep its usual mathematical meaning: $\pi\in\R$ is the constant ratio
between the circumference and the diameter of a circle. With these
conventions, we can omit the ``bar'' which distinguishes a relation or
function from its virtual extension, for its argument indicates which is the
case. Since we are considering $\R\sc\RR$, we will also write $0\in\RR$,
$1\in\RR$, $-3\in\RR$ or~$e\in\RR$ (Neper's number), instead of $\EXT0$,
$\EXT1$, $\EXT{-3}$ or~$\EXT e$. Furthermore, all syntactic and notational
conventions universally accepted for real numbers will be maintained
in~$\RR$. For example:
$$
    \a\b^3 + \ga = [\a(\b^3)] + \ga,
$$
$$
    \sin^2\x = (\sin\x)^2 \qquad\hbox{and}\qquad \sin\x^2 = \sin(\x^2).
$$

The Virtual Extension Theorem (VET, Ref.~1) guarantees the consistence of
identi\-fica\-tions and notations above, provided that we observe relational
symbols which have the underlying logical connectives `not' and~`or'. In this
case, we will maintain the ``bar'' which distinguishes a relation in~$\R$
from its virtual extension in~$\RR$, to avoid ambiguity. For instance, the
expression
$$
    \a\neq\b
$$
means:
$$
    \hbox{$\a=\b$ is false,}
$$
whereas the virtual extension of the relation $\neq\ \sc\R^2$ is represented
by `$\EXT\neq$'. According to item~(ii) of the~VET, the condition
$\a\EXT\neq\b$ is sufficient but not necessary for the validity of
$\a\neq\b$, i.e.:
$$
    \a\EXT\neq\b \IMP \a\neq\b.
$$

Analogously, the expression
$$
    \a\le\b
$$
means
$$
    \a<\b\OR\a=\b,
$$
whereas the virtual extension of the relation $\le\ \sc\R^2$ is represented
by `$\EXT\le$'. By item~(iv) of the~VET, the condition $\a\le\b$ is
sufficient but not necessary for the validity of $\a\EXT\le\b$, i.e.:
$$
    \a\le\b \IMP \a\EXT\le\b.
$$

To avoid ambiguities, we will not write
$$
    \a\EXT{\not\le}\b,
$$
or any other combination of symbols which at the same time involves the bar
that indicates virtual extension and the logical connectives `not' and~`or'.
However, there is no risk of misunderstandings when writing:
$$
    \a\not<\b,
$$
which simply means
$$
    \hbox{$\a<\b$ is false},
$$
or while combining many simultaneous relations on the same variable:
$$
    \a\ne\b\le\ga,
$$
which means simply
$$
    \a\ne\b \AND \b\le\ga.
$$

The multiplicative inverse of a virtual number $\a\in\RR$ will be indicated
by~$\a^{-1}$ or
$$
    {1\over\a}.
$$
It is important to note that those symbols are meaningful just when
$\a\EXT\neq0$, since only virtuals satisfying this condition admit
multiplicative inverse. Otherwise, the fraction above will be treated exactly
like
$$
    {1\over0}.
$$
By the VET, we have that if $\a>0$ or $\a<0$ then $\a\EXT\neq0$, so $\a$ is
inversible.

We will introduce now a special notation for some ``non-real virtuals'':
the class of the sequence $(1,2,3,\ldots)$ will be represented by
$\8\in\RR$, and the multiplicative inverse of this virtual number will be
represented by~$\d\in\RR$, i.e.:
$$
    \d = \<1,{1\over2},{1\over3},\ldots>\in\RR.
$$
Besides, for any $\a=\<a_1,a_2,a_3,\ldots>\in\RR$, we define:
$$
    \pm\a = \<-a_1,+a_2,-a_3,\ldots>\in\RR
$$
and
$$
    \mp\a = \<+a_1,-a_2,+a_3,\ldots>\in\RR.
$$
Those definitions are clearly independent of the representative sequence of
the class $\a\in\RR$.

We will often use the following ternary relation: a real number $b$ {\gr is
between\/} two other reals $a$ and~$c$ when $a\le b\le c$ or~$c\le b\le a$.
This ternary relation defined on~$\R$ extends, as any other $n$-ary relation,
to the set of virtual numbers~$\RR$. We will say $\b$ is between $\a$
and~$\ga$ when the triple $(\a,\b,\ga)$ satisfies that extension. Thus, if
$\a\EXT\le\b\EXT\le\ga$ or~$\ga\EXT\le\b\EXT\le\a$ then $\b$ is between $\a$
and~$\ga$, but the reciprocal statement does not hold: zero is between $\pm1$
and~$\mp1$, despite $\pm1\EXT\le0\EXT\le\mp1$ and $\mp1\EXT\le0\EXT\le\pm1$
being both false.

The aim of all conventions and notations stated in this section is to make
the manip\-ula\-tion of virtual numbers as simple and intuitive as possible.
Those conventions, along with the~VET, allow us to work with virtual numbers
``as if they were real'', without any risk of ambiguities or contradictions.
For instance, the following statements in~$\RR$ immediately come from
the~VET:
$$
    \l > 0 \IMP \m + \l > \m,
$$
$$
    \x > 0 \IMP -\x < 0,
$$
$$
    (\a+\b)^2=\a^2+2\a\b+\b^2,
$$
and
$$
    \sin^2\w + \cos^2\w = 1.
$$
Therefore:
$$
    \8+1 > \8, \qquad 2-\d < 2, \qquad (\8+\d)^2 = \8^2+2+\d^2,
$$
and
$$
    \sin^2\8 + \cos^2\8 = 1.
$$
As another example of application of the~VET, we also have:
$$
    \hbox{$\sin(\n\pi) = 0$, for every virtual integer $\n\in\ZZ$}.
$$
Hence:
$$
    \sin(\8\pi) = 0,
$$
for $\8\in\ZZ\sc\RR$.

\tit{III. Absolute Finitude}

We will say that a virtual number is {\gr infinitesimal} when its absolute
value is less than each real positive number. In other words, $\ep\in\RR$ is
infinitesimal when $|\ep|<x$, for every real number $x>0$.

As an illustration, the virtuals $\d$, $\pm\d$ and~$\mp\d$ are all
infinitesimal. We also have that $\sin\d$ is infinitesimal, since
$|\sin\d\,|<|\d|$. It is clear that we obtain an equivalent definition
substituting ``$|\ep|<x$'' by ``$|\ep|\le x$'' or ``$|\ep|\EXT\le x$'' in the
condition above.

A virtual number will be called {\gr finite\/} when its absolute value is
less than some real number, i.e., $\l\in\RR$ is finite if there exists a real
$x$ such that $|\l|<x$.

For instance, the virtual $\pm1$ is finite, since $|\pm1|=1<2$. The virtuals
$\cos\a$ and~$\sin\a$ are also finite, for every~$\a\in\RR$, since the~VET
guarantees that:
$$
    |\cos\a\,|<2 \qquad\hbox{and}\qquad |\sin\a\,|<2.
$$
In particular, we have that $\cos\d$ and~$\sin\8$ are both finite.

According to the definitions above, every infinitesimal is finite, but the
reciprocal statement does not hold, as the following fact shows:

{\gr The unique infinitesimal real number is the zero, but all real numbers
are finite.}

\PF For every non-zero real number~$x$, it holds that:
$$
    |x| > {|x|\over 2},
$$
hence $x\neq 0$ is not infinitesimal.

On the other hand, for any real $y$, we have:
$$
    |y| < |y| + 1,
$$
so $y$ is finite.\FP

We will say that a virtual is {\gr infinite\/} when it is not finite, i.e.,
$\w\in\RR$ is infinite when $|\w|\not<x$, for every real~$x$. Moreover, we
will say that $\w>\R$ when $\w>x$ for each $x\in\R$; and that $\w<\R$ when
$\w<x$ for each~$x\in\R$.

With those definitions, if $\w>\R$ or $\w<\R$ then $\w$ is infinite, but
there exists infinite virtuals which are neither more nor less than~$\R$. In
addition, if `$\8$' has been read as ``infinity'', then we should clearly
distinguish this noun which denotes one particular virtual number from the
adjective ``infinite'', which merely means ``not finite'', and that can be
applied to many virtual numbers different from~$\8\in\RR$.

Exemplifying: the virtuals $\8$, $-\8$, $\pm\8$ and~$\mp\8$ are all infinite.
Besides, we have $\8>\R$ and $-\8<\R$, but $\pm\8\not>\R$ and $\pm\8\not<\R$.

We present below some basic relations between finitude and the operations of
addition and multiplication, whose proofs can be compared to the traditional
techniques to formal manipulation of limits.

{\gr The sum of two infinitesimals is also infinitesimal.}

\PF We have $|a+b| \le |a| + |b|$, for any $a,b\in\R$. Thus the VET
guarantees that:
$$
    |\ep+\dt|\EXT\le|\ep|+|\dt|,
$$
for any virtuals $\ep$ and~$\dt$. If $\ep$ and~$\dt$ are infinitesimals then
$|\ep|<x/2$ and~$|\dt|<x/2$, for every real $x>0$. Hence, $|\ep+\dt|$ is less
than each positive real~$x$.\FP

{\gr The product of an infinitesimal and a finite virtual is infinitesimal.}

\PF By the VET, he have $|\ep\l|=|\ep||\l|$, for any virtuals $\ep$ and~$\l$.
If $\l$ is finite then there exists a positive real~$y$ such that $|\l|<y$.
If $\ep$ is infinitesimal then $|\ep|<x/y$, for each $x>0$. Therefore
$|\ep\l|$ is less than every positive real~$x$.\FP

This result shows that {\gr the inverse of an infinitesimal, if it exists, is
necessarily infinite}. Otherwise, we would conclude that $1\in\R$ is
infinitesimal.

{\gr The sum and the product of two finite virtuals are both finite.}

\PF If $\l$ and~$\m$ are finite then there exist reals $x$ and~$y$ such that
$|\l|<x$ and $|\m|<y$. Thus we have:
$$
    |\l+\m|<x+y   \qquad\hbox{and}\qquad   |\l\m|<xy,
$$
hence $\l+\m$ and $\l\m$ are both finite.\FP

As a consequence, we have that {\gr the sum of an infinite virtual $\w$ and a
finite virtual $\l$ is necessarily infinite}, for if $\w+\l$ were finite then
$\w=(\w+\l)-\l$ would also be.

However, it is not true that the product of an infinite virtual and a finite
virtual is necessarily infinite, since the finite factor might be
infinitesimal. Nevertheless, we have:

{\gr The product of an infinite virtual and a non-zero real is also
infinite.}

\PF Let $\w$ be an infinite virtual and $a$ a real number different from
zero. If $a\w$ were finite then $|a\w|<x$, for some real~$x$. In this case,
we would have $|\w|<x/|a|$, and so $\w$ would not be infinite.\FP

{\gr If $\w>\R$ or $\w<\R$ then $\w$ is inversible, and its inverse is
infinitesimal.}

\PF If $\w>\R$ then $\w>0$, hence $\w$ is inversible. Furthermore, we have
that $\w>(1/x)$, for every positive~$x$, so
$$
    0<{1\over\w}<x,
$$
for every positive~$x$. The case $\w<\R$ is analogous.\FP

Many other results involving infinitely big virtuals are easily obtained. For
instance:

{\gr If $\w>\R$ and $\ps>\R$ then $\w+\ps>\R$ and $\w\ps>\R$.}

{\gr If $\w<\R$ and $\ps<\R$ then $\w+\ps<\R$ and $\w\ps>\R$.}

{\gr If $\w>\R$ and $\ps<\R$ then $\w\ps<\R$.}

\tit{IV. Proximity}

We will now introduce an equivalence relation on the set~$\RR$ of virtual
numbers: we will say that $\a$ {\gr is infinitely close to~$\b$}, or simply
that $\a$ {\gr is near~$\b$}, when the difference between them is
infinitesimal. In addition, we will write $\a\=\b$ to indicate that
$\a\in\RR$ is near $\b\in\RR$. In other words, $\a\=\b$ when $\a-\b$ is
infinitesimal.

It is easily seen that $\=$ is a symmetric and reflexive relation. To verify
its transitivity, it is enough to note that:
$$
    \a - \ga = (\a - \b) + (\b - \ga),
$$
hence $\a\=\b$ and $\b\=\g$ implies $\a\=\g$, since the sum of two
infinitesimals is also infinitesimal. So we have that {\gr proximity\/} is in
fact an equivalence relation on~$\RR$.

According to the definition above, infinitesimals are exactly the virtuals
infinitely close to zero. The following proposition shows that the proximity
relation preserves the finitude attributes:

{\gr If $\a\=\b$ then:

\(i) $\a$ is infinitesimal \IFF\ $\b$ is infinitesimal\/{\rm;}

\(ii) $\a$ is finite \IFF\ $\b$ is finite\/{\rm;}

\(iii) $\a$ is infinite \IFF\ $\b$ is infinite.}

\PF The two first statements follow from $ \a = (\a - \b) + \b$, since the
sum of two infinitesimals is also infinitesimal, and the sum of two finites
is also finite. The third statement is an immediate consequence of the
second.\FP

Two distinct real numbers cannot be near each other, since the difference
between them, being a non-zero real, cannot be infinitesimal. In other words:
$$
\hbox{\gr if $x\=y$ then $x=y$.}
$$
As a consequence, we have that:
$$
\hbox{\gr if $x\=\a\=y$ then $x=y$.}
$$
Which means: none of the virtuals can be infinitely close to two distinct
reals at the same time.

Not every virtual is near some real, i.e., there exists $\x\in\RR$ such that
$\x\not\=x$, for each $x\in\R$. For example, $\8$ and~$\pm1$ are not
infinitely close to any real number.

It is also not true that $\w>\R$ and $\ps>\R$ implies $\w\=\ps$. To
illustrate, $\8>\R$ and $\8^2>\R$, but $\8\not\=\8^2$, since
$\8^2-\8=\8(\8-1)>\R$ is not infinitesimal.

The following proposition corresponds to the ``Confront Theorem'' of the
traditional formalization:

{\gr If $\b$ is between $\a$ and~$\ga$ and $\a\=\ga$ then $\a\=\b\=\ga$.}

\PF In~$\R$ we have that if $b$ is between $a$ and~$c$, then:
$$
    |b-a| \le |c-a|.
$$
Thus, the VET guarantees that if $\b$ is between $\a$ and~$\ga$ then:
$$
    |\b-\a| \EXT\le |\ga-\a|.
$$
Therefore $\ga\=\a$ implies $\b\=\a$.\FP

Addition preserves the proximity between virtual numbers:

{\gr If $\a_1\=\a_2$ and $\b_1\=\b_2$ then $\a_1 + \b_1 \= \a_2 + \b_2$.}

\PF It is enough to note that:
$$
    (\a_1 + \b_1) - (\a_2 + \b_2) = (\a_1 - \a_2) + (\b_1 - \b_2),
$$
and remember that the sum of two infinitesimals is also infinitesimal.\FP

In particular, we have that $\a\=x$ and $\b\=y$ implies $\a+\b\=x+y$. This
corresponds to the statement ``The limit of the sum equals the sum of the
limits'' from traditional Calculus.

About the relation between multiplication and proximity, we have:

{\gr If $\l_1\=\l_2$ and $\m_1\=\m_2$, with $\l_1$ and~$\m_1$ finite, then
$\l_1\m_1 \= \l_2\m_2$.}

\PF For every quadruple of virtuals $\l_1$, $\m_1$, $\l_2$ and $\m_2$, it
holds that:
$$
    \l_1\m_1 - \l_2\m_2 = \l_1\m_1-\l_1\m_2+\l_1\m_2-\l_2\m_2 =
      \l_1(\m_1-\m_2) + (\l_1-\l_2)\m_2.
$$
Since $\l_1$ and~$\m_2$ are finite, and $(\m_1-\m_2)$ and~$(\l_1-\l_2)$ are
infinitesimals, we conclude that $(\l_1\m_1 - \l_2\m_2)$ is also
infinitesimal.\FP

A counterexample having one of the factors infinite: let $\l_1=1$,
$\l_2=1+\d$, and $\m_1=\m_2=\8$. In this case:
$$
    \l_1\m_1 - \l_2\m_2 = \8 - (1+\d)\8 = -\d\8 = -1,
$$
so $\l_1\m_1 \not\= \l_2\m_2$.

Particularizing the result above, we obtain that $\a\=x$ and $\b\=y$ implies
$\a\b\=xy$. (The limit of the product equals the product of the limits.)

For multiplicative inversibility we have:

{\gr If $\a\=x\ne0$ then $\a$ is inversible and $1/\a\=1/x$.}

\PF If $x>0$ then $x/2<\a<2x$, and if $x<0$ then $2x<\a<x/2$. In both cases
$\a$ is inversible with $1/\a$ finite. Since
$$
    {1\over\a} - {1\over x} = {1\over \a x}(x-\a),
$$
we have $1/\a\=1/x$, for $1/(\a x)=(1/\a)(1/x)$ is finite and $(x-\a)$ is
infinitesimal.\FP

Therefore $\a\=x\neq0$ and~$\b\=y$ implies:
$$
    {\b\over\a} \= {y\over x}.
$$

We present below some examples that can be compared to the corresponding
``computations'' in the traditional formalization:
$$
    { (5 + \d)^2 - 25 \over \d } =
    { 25 + 10\d + \d^2 - 25 \over \d } =
    { 10\d + \d^2 \over \d } = 10 + \d \= 10.
$$
$$
    { 2\8^3 + 4\8^2 -1 \over \8^3 - 5 } =
    { \8^3\left( 2 + {4\over\8} - {1\over\8^3}\right) \over
      \8^3\left( 1 - {5\over\8^3}\right) } =
    { 2 + 4\d - \d^3 \over 1 - 5\d^3  } \= 2.
$$
$$\eqalign{\8 - \sqrt{\8^2+1} &=\left(\8 - \sqrt{\8^2+1}\right)
    {\left(\8 + \sqrt{\8^2+1}\right)\over\left(\8 + \sqrt{\8^2+1}\right)}\cr
  &={\8^2 - \left(\8^2+1\right)\over\left(\8 + \sqrt{\8^2+1}\right)}\cr
  &={ -1 \over \left(\8 + \sqrt{\8^2+1}\right)}\cr
  &\= 0,}
$$
hence $\8\=\sqrt{\8^2+1}$.

\tit{V. Continuity}

A function $\f\:\D\to\C$ with $\D\sc\R$ and~$\C\sc\R$ is continuous when its
virtual extension:
$$
    \ff\:\DD\to\CC
$$
preserves the proximity relation on~$\RR$. As established in Sec.~II, we will
omit the bar which distinguishes $\f$ from its virtual extension, in order to
simplify notation.

We will say that $\f$ {\gr is continuous at\/} $x\in\D$ when
$$
    \V \a\in\DD,\ \a\=x \IMP \f(\a)\=\f(x).
$$
This condition is equivalent to the usual definition of continuity of~$\f$ at
a point of its domain. We will say that {\gr $\f$ is continuous\/} when it is
continuous at each~$x\in\D$. This definition is equivalent to the
``pointwise'' continuity of traditional Calculus.

{\gr If $\f$ and~$\g$ are real functions defined on the same domain $\D\sc\R$
and continuous at~$x\in\D$ then $\f+\g$ is also continuous at~$x$.}

\PF Since $(\f+\g)(a)=\f(a)+\g(a)$, for every $a\in\D$, The VET guarantees
that $(\f+\g)(\a)=\f(\a)+\g(\a)$, for every $\a\in\DD$. Now, it is enough to
remember that $\f(\a)\=\f(x)$ and~$\g(\a)\=\g(x)$ implies
$\f(\a)+\g(\a)\=\f(x)+\g(x)$.\FP

Therefore, {\gr the sum of two continuous functions defined on the same
domain is also continuous}.

Analogously one verifies that:

{\gr If $\f$ and~$\g$ are real functions defined on the same domain $\D\sc\R$
and continuous at~$x\in\D$ then $\f\g$ is also continuous at~$x$. Besides, if
$\g(x)\neq0$ then $\f/g$ is also continuous at~$x$.}

Consequently, {\gr the product and the quotient of two continuous functions
defined on the same domain are both continuous}.

The intuitive argument for the following result and its formal demonstration
are ident\-ical:

{\gr If two real functions $\f$ and~$\g$ are continuous and make a chain then
the composite function $(\g\o\f)$ is also continuous.}

\PF If $x\in\D$ and $\a\in\DD$ with $\a\=x$ then $\f(\a)\=\f(x)$, for $\f$ is
continuous. Since $\g$ is also continuous, we have $\g[\f(\a)]\=\g[\f(x)]$.
Therefore the composite $(\g\o\f)$ is continuous.\FP

Using virtual numbers, we can state a more restrictive condition of
continuity than the one used so far: we say that $\f\:\D\to\C$ is {\gr
uniformly continuous\/} on $\A\sc\D$ when
$$
    \V \a,\b\in\AA,\ \a\=\b \IMP \f(\a)\=\f(\b).
$$
This condition is equivalent to the uniform continuity of traditional
Mathematical Analysis. We can demonstrate this fact by the usual methods,
keeping in mind that virtual numbers are classes of sequences of real
numbers. However, we will present a proof based on the~VET, which illustrates
how it can substitute those traditional techniques.

{\gr If $\f:\D\to\C$ and~$\A\sc\D$, then the following conditions are
equivalent:

\(i) $\V \a,\b\in\AA,\ \a\=\b \IMP \f(\a)\=\f(\b)$.

\(ii) $\V r>0,\E s>0,\V a,b\in\A,\ |a-b|<s \IMP |\f(a)-\f(b)|<r$.}

\PF If condition (ii) holds, then, by the~VET, we have:
$$
    \V r>0,\E s>0,\V \a,\b\in\AA,\ |\a-\b|<s \IMP |\f(\a)-\f(\b)|<r.
$$
Since $\a\=\b$ implies $|\a-\b|<s$, for every $s>0$, we conclude that
$\a\=\b$ implies
$$
    \V r>0,\ |\f(\a)-\f(\b)|<r.
$$
That means, $\a\=\b$ implies $\f(\a)\=\f(\b)$. Therefore (i) follows
from~(ii).

On the other hand, if condition (ii) does not hold, then:
$$
    \E r>0,\V s>0,\E a,b\in\A,\ |a-b|<s \AND |\f(a)-\f(b)| > r.
$$
By the VET, we obtain:
$$
    \E r>0,\V\s>0,\E\a,\b\in\AA,\ |\a-\b|<\s \AND |\f(\a)-\f(\b)| > r.
$$
Now taking $\s\=0$ ($\d$ for example), we get:
$$
    \E r>0,\E\a,\b\in\AA,\ \a\=\b  \AND |\f(\a)-\f(\b)| > r,
$$
i.e., there exist $\a,\b\in\AA$ with $\a\=\b$ and $\f(\a)\not\=\f(\b)$, which
is the negation of~(i).\FP

We can verify the uniform continuity of a real function directly by the
proposed definition, i.e., by condition (i) above. For example, the function
cosine is uniformly continuous on~$\R$. To prove this, we note that:
$$
    |\cos x - \cos y\,| \ \le\ |x - y|,
$$
for any reals $x$ and~$y$. Therefore, by the~VET, we have
$$
    |\cos\x - \cos\y\,| \ \EXT\le\ |\x - \y|,
$$
for any virtuals $\x$ and~$\y$. Hence, $\x\=\y$ implies $\cos\x\=\cos\y$.

The notion of pointwise continuity and that of uniform continuity can both be
defined using only the proximity relation between virtual numbers. Many other
topological notions depend just on that relation, which motivates the
following definition:

A {\gr virtual topology\/} on a generic set~$\A$ is an equivalence relation
$\=$ on its virtual extension~$\AA$. We will say that $\=$ is {\gr
separable\/} when there exists at the most one element of $\A\sc\AA$ in each
equivalence class, i.e., when $x\=y$ with $x,y\in\A$ implies $x=y$. As we
have seen, the proximity between virtual numbers defined in the previous
section is a separable virtual topology on~$\R$.

The notions of {\gr metric space\/} and {\gr topological space\/} are
generalizations of the traditional definition of continuity. A virtual
topology allows us to introduce topological notions on any set in a simpler,
more direct and intuitive way. In addition, we can define uniform continuity
using just a virtual topology (as we have seen above), which cannot be done
by using just the traditional topological structure.

We define below two familiar topological concepts on any {\gr virtual
topological space}. They will be used further for the particular case of real
line.

We will say that $x\in\A$ is an {\gr interior point\/} of a set $\B\sc\A$
when $\a\=x$ implies~$\a\in\BB$, and we will say that $\B\sc\A$ is an {\gr
open set\/} when each $x\in\B$ is an interior point of~$\B$.
\vfill\eject
\tit{VI. Derivation}

To simplify notation, we will now define another relation between virtual
numbers: we will say that $\a$ and~$\b$ are {\gr neighbours\/} when $\a\=\b$
and $\a\EXT\neq\b$. Moreover, we will write $\a\@\b$ to indicate that
$\a\in\RR$ and $\b\in\RR$ are neighbours. This neighbourliness relation is
symmetric, but it is neither reflexive nor associative, and so it is not an
equivalence relation on~$\RR$. According to this definition, the inversible
infinitesimals are exactly the neighbours of zero.

The following proposition illustrates the use of the neighbourliness
relation:

{\gr If $\ep\in\RR$ is an invesible infinitesimal, i.e., if $\ep\@0$, then:}
$$
    {\sin\ep\over\ep}\=1.
$$

\PF It holds in~$\R$ that:
$$
    0 \neq |x| < {\pi\over2} \IMP \cos x < {\sin x\over x} < 1.
$$
Therefore, by the VET, we get in~$\RR$:
$$
    0 \EXT\neq |\x| < {\pi\over2} \IMP \cos\x < {\sin\x\over\x} < 1.
$$
Thus we have, for any $\ep\@0$:
$$
    \cos\ep < {\sin\ep\over\ep} < 1,
$$
so $(\sin\ep)/\ep$ is between $\cos\ep$ and~1. Since the cosine function is
continuous, we conclude that $\cos\ep\=\cos0=1$, hence:
$$
    {\sin\ep\over\ep}\=1.\FP
$$

Let now $\f\:\D\to\C$ be a function with $\D\sc\R$ and~$\C\sc\R$. We will say
that $\f$ is {\gr derivable at~$x\in\D$} when there exists $m\in\R$ such
that
$$
    \V\a\in\DD,\ \a\@ x \IMP { \f(\a) - \f(x) \over \a - x } \= m.
$$
That condition is equivalent to the traditional definition of derivability at
$x\in\D$, not presuming that $x$ is an interior point of the domain of~$\f$.
We will say simply that {\gr $\f$ is derivable\/} when it is derivable at
each point of its domain.

If $\f\:\D\to\C$ is a derivable function then there exists exactly one
$m_x\in\R$ which satisfies the derivability condition for each $x\in\D$. Thus
we define the {\gr derivative} $\f'$ of~$\f$ by:
$$\eqalign{
    \f'\:&\D\to\R\cr
         &x\mapsto m_x.\cr
}$$

Let $y$ be the real variable depending on~$x$ by a function~$\f$, i.e.:
$$
    y = \f(x).
$$
If $x$ is an interior point of the domain of~$\f$, we can consider an {\gr
infinitesimal variation\/} of~$x$ as a virtual number~$dx\@0$. The
corresponding variation of~$y$ is:
$$
    dy = \f(x+dx) - \f(x).
$$
So, it is clear that $\f$ is derivable at~$x$ when, for any virtual
variation~$dx\@0$, the quotient of~$dy$ and~$dx$ is infinitely close to the
same real number. In this case, we have:
$$
    {dy\over dx} \= \f'(x),
$$
which can be seen as a formalization of the ``manipulations of infinitesimal
variations'' presented in the classical notation due to Leibniz.

Some examples:

(i) If $\f(x)=k$, for every $x\in\R$, then:
$$
    {dy\over dx} = {\f(x+dx)-\f(x)\over dx} = {k - k \over dx} = 0,
$$
hence, any constant function is derivable, and its derivative is the constant
zero function.

(ii) If $\f$ is the identity function, i.e., if $y=x$, then:
$$
    {dy\over dx} = {(x+dx) - x \over dx} = {dx\over dx} = 1,
$$
so the identity function is derivable and $\f'(x)=1$, for every $x\in\R$.

(iii) Let $y=\sin x$. By the~VET, we have:
$$
    \sin\a-\sin\b = 2\sin{\a-\b\over2}\cos{\a+\b\over2},
$$
for any virtuals $\a$ and~$\b$. Hence:
$$
  \eqalign{{dy \over dx} &={\sin(x+dx)-\sin x \over dx}\cr
    &={2\sin{\displaystyle dx\over\displaystyle 2}\cos{\displaystyle
      2x+dx\over\displaystyle 2}\over dx}\cr
    &={\sin{\displaystyle dx\over\displaystyle 2}\over{\displaystyle
      dx\over\displaystyle 2}}\cos\left(x+{dx\over2}\right)\cr
    &\=\cos x.\cr
}$$
Therefore the sine function is derivable and its derivative is the cosine
function.

We also have the usual {\gr derivation rules:}

{\gr If $\f$ and~$\g$ are real functions derivable at an interior point~$x$
of both domains, then $(\f+\g)$ is also derivable at~$x$, and:}
$$
    (\f+\g)'(x) = \f'(x) + \g'(x).
$$

\PF For every $dx\@0$ we have:
$$
  \eqalign{{(\f+\g)(x+dx)-(\f+\g)(x)\over dx}
    &={\f(x+dx)+\g(x+dx)-\f(x)-\g(x)\over dx}\cr
    &={\f(x+dx)-\f(x)\over dx}+{\g(x+dx)-\g(x)\over dx}\cr
    &\= \f'(x)+\g'(x).\FP\cr}
$$

{\gr If $\f$ and~$\g$ are real functions derivable at an interior point~$x$
of both domains, then $(\f\g)$ is also derivable at~$x$, and:}
$$
    (\f\g)'(x) = \f'(x)\g(x) + \f(x)\g'(x).
$$

\PF For every $dx\@0$ we have:
$$\eqalign{&{(\f\g)(x+dx)-(\f\g)(x)\over dx}=\cr
 &\qquad={\f(x+dx)\g(x+dx)-\f(x)\g(x+dx)+\f(x)\g(x+dx)-\f(x)\g(x)\over dx}\cr
 &\qquad={\f(x+dx)-\f(x)\over dx}\g(x+dx) + \f(x){\g(x+dx)-\g(x)\over dx}\cr
 &\qquad\=\f'(x)\g(x)+\f(x)\g'(x).\FP\cr
}$$

{\gr If $\g\:\D\to\R$ is a real function derivable at an interior point~$x$
of its domain and $\g(x)\ne0$ then $(1/\g)$ is derivable at~$x\in\D$, and:}
$$
    \left({1 \over \g}\right)'\!(x) = - {\g'(x) \over [\g(x)]^2}.
$$

\PF Let $dx\@0$. We have that $\g(x+dx)$ is inversible, since $\g$ continuous
at~$x$ implies $\g(x+dx)\=\g(x)\ne0$. Furthermore:
$$
    { {\displaystyle 1\over\displaystyle\g(x+dx)}
    - {\displaystyle 1\over\displaystyle\g(x)} \over dx}
    = - {\g(x+dx) - \g(x) \over dx} {1 \over \g(x+dx)\g(x)}
    \= - {\g'(x) \over [\g(x)]^2}.\FP
$$

There is a more restrictive condition of derivability than the one used so
far: we will say that $f\:\D\to\R$ is {\gr differentiable at~$x\in\D$} when
there exists $m\in\R$ such that:
$$
    \V \a,\b\in\DD,\ \a\@\b\=x \IMP {\f(\a)-\f(\b) \over \a-\b} \= m.
$$
It is clear that differentiability implies derivability (as defined in this
section). The following theorem explains the relation between these two
concepts:

{\gr Let $\f\:\D\to\R$ be a derivable function defined on an open set
$\D\sc\R$. Its derivative $\f'\:\D\to\R$ is continuous at~$x\in\D$ \IFF\ $\f$
is differentiable at this point~$x$.}

To prove this theorem, we will use the following lemma:

{\gr If $\f\:\D\to\C$ is a derivable function on an open set $\D\sc\R$ then,
for every $\a\in\DD$, there exists $\b\in\DD$ with $\b\@\a$, and}
$$
    { \f(\a) - \f(\b) \over \a - \b } \= \f'(\a).
$$

\PF We will show initially that if $\f$ is derivable on~$\D$ then:
$$
    \V a\in\D,\V r>0,\E b\neq a,\ |b-a|<r \AND
    \left| {\f(a) - \f(b) \over a - b} - \f'(a) \right| <r.
$$
If this statement were false we would have:
$$
    \E a\in\D,\E r>0,\V b\neq a,\ |b-a|<r \IMP
    \left| {\f(a) - \f(b) \over a - b} - \f'(a) \right| > r.
$$
Thus, we would conclude, by the~VET, that:
$$
    \E a\in\D,\E r>0,\V\b\EXT\neq a,\ |\b-a|<r \IMP
    \left| {\f(a) - \f(\b) \over a - \b} - \f'(a) \right| > r,
$$
and so:
$$
    \E a\in\D,\E r>0,\V\b\@ a,\
    \left| {\f(a) - \f(\b) \over a - \b} - \f'(a) \right| > r,
$$
which means $\f$ would not be derivable at~$a\in\D$.

Applying the VET to the statement proved above we get:
$$
    \V\a\in\DD,\V\r>0,\E\b\EXT\neq\a,\ |\b-\a|<\r \AND
    \left| {\f(\a) - \f(\b) \over \a - \b} - \f'(\a) \right| < \r.
$$
Now taking $\r$ infinitesimal ($\d$ for example), we conclude that, for every
$\a\in\DD$, there exists $\b\@\a$ such that:
$$
    {\f(\a) - \f(\b) \over \a - \b} \= \f'(\a).\FP
$$

To prove the theorem, let $\f\:\D\to\C$ be a function derivable on an open
set $\D\sc\R$. If $\f'$ is not continuous at~$x\in\D$, then there exists
$\a\=x$ such that $\f'(\a)\not\=\f'(x)$. Thus, by the lemma above, we
conclude that there exists $\b\@\a\=x$ such that
$$
    {\f(\a) - \f(\b) \over \a - \b} \= \f'(\a) \not\= \f'(x),
$$
which reads: the function $\f$ is not differentiable at~$x$.

On the other hand, the Mean Value Theorem of traditional Calculus guarantees
that for any pair $a,b\in\D$ with~$a\neq b$, there exists $c$ between $a$
and~$b$ such that:
$$
    {\f(a) - \f(b) \over a - b} \= \f'(c).
$$
Thus, by the~VET, we have that for any pair $\a,\b\in\DD$ with~$a\EXT\neq\b$
there exists $\ga$ between $\a$ and~$\b$ such that:
$$
    {\f(\a) - \f(\b) \over \a - \b} \= \f'(\ga).
$$
If the derivative of~$\f$ is a function continuous at~$x\in\D$ then, for any
pair $\a,\b\in\DD$ with~$\a\@\b\=x$, there exists $\ga\=x$ such that
$$
    {\f(\a) - \f(\b) \over \a - \b} \= \f'(\ga) \= \f'(x).
$$
Therefore $\f$ is differentiable at~$x$.

\tit{VII. Relative Finitude}

Let $\a\in\RR$ be an inversible virtual, i.e., $\a\EXT\neq0$. We will say
that $\b\in\RR$ {\gr is of order~$\a$} when the quotient $\b/\a$ is finite.
In addition, we will represent the set of all virtual numbers of order~$\a$
by $\OO\a\sc\RR$. That means:
$$
    \OO\a=\left\{ \b\in\RR\ |\ {\b\over\a} \hbox{\ is finite\ }\right\}.
$$
For instance, $\sin\ep\in\OO\ep$, for any inversible infinitesimal $\ep\@0$.

We will say that $\ga\in\RR$ is {\gr negligible when compared to~$\a$}, if
the quotient $\ga/\a$ is infinitesimal, i.e., if
$$
    {\ga\over\a} \= 0.
$$
We will indicate that by writing:
$$
    \ga\ll\a.
$$
Exemplifying: the square of any inversible infinitesimal $\ep\@0$ is
negligible when compared to~$\ep$ itself:
$$
    \ep\@0 \IMP \ep^2\ll\ep.
$$

According to those definitions, it is obvious that if $\ga$ is negligible
when compared to~$\a$, then $\ga$ is of order~$\a$:
$$
    \ga\ll\a \IMP \ga\in\OO\a,
$$
and that the reciprocal statement does not hold.

Using results about absolute finitude and proximity (like those presented in
Secs.~III and~IV, for instance) we can establish many relations between those
concepts and the notions of relative finitude introduced above. As an
illustration, it is easily shown that any infinitesimal is negligible when
compared to any non-zero real number, and that any finite virtual is
negligible when compared to~$\8$.

The following proposition states the {\gr Taylor's Formula\/} of a real
function as it is actually used by most scientists. Our notation will be like
this:
$$
    \a = \b + \OO\ga,
$$
which merely means:
$$
    \a - \b  \in \OO\ga.
$$

{\gr If $x$ is an interior point of the domain of a function $\f$ which is
differentiable $n+1$ times ($n\in\N$), then for any $\ep\@ 0$ we have:}
$$
    \f(x+\ep) = \f(x) + \f'(x)\ep +\cdots+ {\f^{(n)}(x)\over n!}\ep^n +
                \OO{\ep^{n+1}}.
$$

\PF Let $\I$ be an open interval inside the domain of~$\f$ such that
$x\in\I$. Taylor's Formula with Lagrange's Remainder guarantees that, for
each $a\in\I$, there exists $b$ between $x$ and~$a$ such that:
$$
    \f(a) = \f(x) + \f'(x)(a-x) +\cdots+ {\f^{(n)}(x)\over n!}(a-x)^n +
            {\f^{(n+1)}(b)\over (n+1)!}(a-x)^{n+1}.
$$
Thus, we conclude, by the~VET, that for each $\a\in\II$ there exists $\b$
between $x$ and~$\a$ such that:
$$
    \f(\a) = \f(x) + \f'(x)(\a-x) +\cdots+ {\f^{(n)}(x)\over n!}(\a-x)^n +
             {\f^{(n+1)}(\b)\over (n+1)!}(\a-x)^{n+1}.
$$

Now taking $\a=x+\ep$, it follows that there exists $\b\=x$ such that:
$$
    \f(x+\ep) = \f(x) + \f'(x)\ep +\cdots+ {\f^{(n)}(x)\over n!}\ep^n +
             {\f^{(n+1)}(\b)\over (n+1)!}\ep^{n+1}.
$$
Since $\f^{(n+1)}$ is continuous, we have
$\f^{(n+1)}(\b)\=\f^{(n+1)}(x)\in\R$. Hence:
$$
    {\f^{(n+1)}(\b)\over (n+1)!}\ep^{n+1}\in\OO{\ep^{n+1}}.\FP
$$

According to the definitions introduced in this section, we can consider the
{\gr First Rule of L'Hospital\/} as a criterion for deciding the relative
magnitude between two infinitesimals $\ep$ and~$\dt$ of the kind:
$$
    \ep=\f(\a)\=0 \qquad\hbox{and}\qquad \dt=\g(\a)\@ 0.
$$
Analogously, the {\gr Second Rule of L'Hospital\/} might be seen as a
criterion to match the magnitude of two infinitely big virtuals, i.e., which
are more than~$\R$ or less than~$\R$.

\tit{VIII. Integration}

Our aim in this section is to illustrate the usefulness of relative finitude
notions introduced in the last section. We will begin with a definition of
integral which uses the language of virtual numbers instead of the idea of
limit. However, it is essentially the construction due to Riemann.

There is another way to define the integral as a sum of infinitely many
infinitesimal terms, using the concept of ``virtual sequence''. We will not
go this way here because it would be premature to introduce that notion now.

Let $\f$ be a function defined on an closed interval with extremes
$a,b\in\R$. We define the following relation~$\P$ between pairs $(m,s)$ of
real numbers:

(i) if $a=b$ then $\P(m,s)$ \IFF\ $m>0$ and~$s=0$;

(ii) if $a<b$ then $\P(m,s)$ \IFF\ there exists an {\gr extended partition\/}
of interval $[a,b]$:
$$
    a = x_0 < z_1 < x_1 < z_2 < x_2 < \cdots < x_n = b,
$$
with {\gr norm\/} less then~$m$:
$$
    m > x_i-x_{i-1} \qquad (i=1,\ldots,n),
$$
whose {\gr Riemann sum\/} is equal to~$s$:
$$
    s = \sum_{i=1}^n \f(z_i)(x_i-x_{i-1});
$$

(iii) if $b<a$ then $\P(m,s)$ \IFF\ there exists an extended partition of the
interval $[b,a]$, with norm less than~$m$, whose Riemann sum equals~$-s$.

We will use the virtual extension $\PP$ of the relation~$\P$ above to define
the integrability of the function~$\f$. According to the convention adopted
in Sec.~II, we will omit the bar which distinguishes $\P$ from its virtual
extension.

We say that $\f$ {\gr is integrable\/} on the interval of extremes $a,b\in\R$
when there exists $s\in\R$ such that $\P(\m,\s)$ and~$\m\=0$ implies $\s\=s$.
Informally: $\f$ is integrable when any extended partition infinitely fine
has the Riemann sum near the same real number. Obviously, if $\f$ is
integrable then there exists a unique real number~$s$ which satisfies that
condition. This number~$s$ will be called {\gr integral of~$\f$ between $a$
and~$b$}.

It is not difficult to see that the integrability condition above and the
corresponding definition of integral are equivalent to those due to Riemann.
Therefore, any continuous function is integrable in this sense.

We can think of the integral as a ``finite sum of an infinite quantity of
infinitesimal elements''. The notation universally adopted for the integral
of~$\f$ between $a$ and~$b$ suggests this idea:
$$
    s = \int_a^b \f(x)\,dx.
$$

It is important to note that this notation is founded on the following fact:
the infinitesimal elements which we ``add'' during the integration process
can be approximated by the product $\f(x)\,dx$, without altering the final
result. This is, essentially, the Fundamental Theorem of Calculus.

To explain this assertion better, let $\f$ be a real continuous function,
defined on an open interval, and $a$ a generic point in this interval. Since
$f$ is integrable, we can define:
$$
    s = \g(x) = \int_a^x \f(t)\,dt.
$$
Thus, for each infinitesimal increment $dx$ of variable~$x$, we have a
corresponding infinitesimal variation of~$s$:
$$
    ds = \g(x+dx) - \g(x).
$$
Those are the elements $ds$ which we ``add'' during the integration process.

The Fundamental Theorem of Calculus states that $\g$ is derivable and its
derivative equals the function~$\f$:
$$
    \g'(x) = \f(x).
$$
This statement is equivalent to
$$
    {ds\over dx} \= \f(x),
$$
which can be rewritten as:
$$
    \left[{ds\over dx} - \f(x)\right] \= 0,
$$
or:
$$
    \left[{ds - \f(x)\,dx \over dx}\right] \= 0.
$$
That means, the Fundamental Theorem of Calculus is equivalent to the
following statement:

{\gr If $\f$ is a continuous function, then, for any $dx\@0$, the error in
the approximation:
$$
    ds \= \f(x)\,dx
$$
is negligible when compared to~$dx$, i.e.:}
$$
    \left[ds - \f(x)\,dx \right] \ll dx.
$$

This assertion admits a quite intuitive graphic interpretation.

Every time we calculate a finite quantity by integration, we use an
approximation for the corresponding infinitesimal element. In order not to
change the final result, it is necessary that the approximation error be
negligible when compared to the infinitesimal element to which we are
approximating. For instance, let us consider the calculation of lengths,
areas, and volumes of some simple geometric objects:

If $\f\:[a,b]\to\R$ is a strictly positive derivable function then:

(i) The area of the plan region under the graph of~$\f$ is given~by:
$$
    A = \int da = \int_a^b\f(x)\,dx.
$$

(ii) The volume of the solid obtained by the revolution of the region under
the graph of~$\f$ around the $x$-axis is given~by:
$$
    V = \int dv = \int_a^b\pi[\f(x)]^2\,dx.
$$

(iii) The length of the graph of~$\f$ is given~by:
$$
    L = \int d\ell = \int_a^b\sqrt{1 + [\f'(x)]^2}\,dx.
$$

(iv) The area of the surface obtained by the revolution of the graph of~$\f$
around the $x$-axis is given~by:
$$
    S = \int ds = \int_a^b 2\pi \f(x)\sqrt{1 + [\f'(x)]^2}\,dx.
$$

The derivative of the function $\f$ is not present in the integrands of the
two first formulae above, but it is indispensable in the last two.

In the first formula, we can approximate the {\gr element of area~$da$} by
$\f(x)\,dx$ because the approximation error is negligible when compared
to~$da$ itself:
$$
    da = \f(x)\,dx + \ep, \qquad\hbox{with}\qquad \ep\ll da.
$$
Analogously, in the second formula, if we approximate:
$$
    dv \= \pi[\f(x)]^2\,dx,
$$
we make a negligible error when compared to the {\gr volume element~$dv$}
itself, which does not alter the final result~$V$.

On the other hand, if we approximate the {\gr element of length~$d\ell$}
simply by~$dx$, the error would not be negligible opposed to~$d\ell$ itself,
and the integral formula thus obtained:
$$
    L = \int_a^b dx = b - a
$$
would not hold. Taking
$$
    d\ell \= \sqrt{1 + [\f'(x)]^2}\,dx,
$$
we commit a negligible error when compared to~$d\ell$, hence it does not
alter the length~$L$ obtained from formula~(iii).

The same way, approximating the area $ds$ of the {\gr surface element\/}
merely by
$$
    2\pi\f(x)\,dx,
$$
the error is not negligible opposed to~$d\ell$ itself, and the corresponding
integral formula:
$$
    S = 2\pi \int_a^b \f(x)\,dx
$$
is wrong. Formula (iv) is right because the error in the approximation:
$$
    ds \= 2\pi\f(x)\,d\ell \= 2\pi\f(x)\,\sqrt{1 + [\f'(x)]^2}\,dx
$$
is negligible when compared to~$ds$ itself.

\tit{Reference}

\item{$^1$}S.~F.~Cortizo: ``Virtual Extensions'', to appear (1995).
\bye